\documentclass[reqno]{amsart} 
\usepackage{amssymb}   
\usepackage{url} 
\usepackage{graphicx}
\def\qmod#1#2{{\hbox{}^{\displaystyle{#1}}}\!\big/\!\hbox{}_{
\displaystyle{#2}}}

\def\resto#1#2{{
#1\hskip 0.4ex\vline_{\hskip 0.2ex\raisebox{-0.2ex}
{{${\scriptstyle #2}$}}}}}

\def\C{{\mathbb C}}

\def\P{{\mathbb P}}

\def\R{{\mathbb R}}

\def\Z{{\mathbb Z}}

\def\qed {\hfill\vrule height6pt width6pt depth0pt \bigskip}

\def\map{\longrightarrow}
\def\textmap#1{\mathop{\vbox{\ialign{
                                 ##\crcr
     ${\scriptstyle\hfil\;\;#1\;\;\hfil}$\crcr
     \noalign{\kern 1pt\nointerlineskip}
     \rightarrowfill\crcr}}\;}}
\newcommand{\cal}{\mathcal}
\def\textlmap#1{\mathop{\vbox{\ialign{
                                 ##\crcr
     ${\scriptstyle\hfil\;\;#1\;\;\hfil}$\crcr
     \noalign{\kern-1pt\nointerlineskip}
     \leftarrowfill\crcr}}\;}}

\def\textlmap#1{\mathop{\vbox{\ialign{
                                  ##\crcr
      ${\scriptstyle\hfil\;\;#1\;\;\hfil}$\crcr
      \noalign{\kern-1pt\nointerlineskip}
      \leftarrowfill\crcr}}\;}}

\def\g{{\mathfrak g}}

\def\Rg{{\mathfrak R}}

\newtheorem{sz}{Satz}[section]
\newtheorem{thry}[sz]{Theorem}
\newtheorem{pr}[sz]{Proposition}
\newtheorem{re}[sz]{Remark}

\begin{document}
\def\Pr{{\rm Pr}}
\def\ASD{{\rm ASD}}
\def\tr{{\rm Tr}}
\def\End{{\rm End}}
\def\Aut{{\rm Aut}}
\def\Spin{{\rm Spin}}
\def\U{{\rm U}}
\def\SU{{\rm SU}}
\def\SO{{\rm SO}}
\def\PU{{\rm PU}}
\def\GL{{\rm GL}}
\def\spin{{\rm spin}}
\def\su{{\rm su}}
\def\so{{\rm so}}
\def\ub{\underbar}
\def\pu{{\rm pu}}
\def\Pic{{\rm Pic}}
\def\Iso{{\rm Iso}}
\def\NS{{\rm NS}}
\def\Hom{{\rm Hom}}
\def\Aut{{\rm Aut}}
\def\h{{\germ h}}
\def\Herm{{\rm Herm}}
\def\Vol{{\rm Vol}}
\def\pf{{\bf Proof: }}
\def\id{{\rm id}}
\def\deg{{\rm deg}}
\def\i{{\germ i}}
\def\im{{\rm im}}
\def\rk{{\rm rk}}
\def\ad{{\rm ad}}
\def\h{{\bf H}}
\def\coker{{\rm coker}}
\def\dbar{\bar{\partial}}
\def\Lo{{\Lambda_g}}
\def\niq{=\kern-.18cm /\kern.08cm}
\def\Ad{{\rm Ad}}
\def\RSU{\R SU}
\def\ad{{\rm ad}}
\def\dva{\bar\partial_A}
\def\da{\partial_A}
\def\p{\mathrm p}
\def\sp{\Sigma^{+}}
\def\sm{\Sigma^{-}}
\def\spm{\Sigma^{\pm}}
\def\smp{\Sigma^{\mp}}
\def\Dr{{\raisebox{0.15ex}{$\not$}}{\hskip -1pt {D}}}
\def\Tors{{\rm Tors}}
\def\st{{\rm st}}
\def\s{{\rm s}}
\def\oo{{\scriptstyle{\mathcal O}}}
\def\ooo{{\scriptscriptstyle{\mathcal O}}}
\def\sw{Seiberg-Witten }
\def\pa{\partial_A\bar\partial_A}
\def\gr{{\scriptscriptstyle|}\hskip -4pt{\g}}
\def\subsetint{{\  {\subset}\hskip -2.45mm{\raisebox{.28ex}
{$\scriptscriptstyle\subset$}}\ }}
\def\ra{\rightarrow}
\def\pst{{\rm pst}}
\def\sst{{\rm sst}}
\def\td{{\rm td}}
\def\kod{{\rm kod}}
\def\degmax{{\rm degmax}}
\def\red{{\rm red}}
\def\reg{{\rm reg}}

\title[Gauge theoretical methods]{Gauge theoretical methods in the classification of non-K\"ahlerian surfaces}
\author{Andrei Teleman}
\address{LATP, CMI, Universit\'e de Provence, 39 Rue Fr\'ed\'eric Joliot-Curie, 13453 Marseille cedex 13, France} 
\email{teleman@cmi.univ-mrs.fr} 
\date{\today}
\begin{abstract}
 The classification of  class VII surfaces is a  very difficult classical problem in complex geometry. It is considered by experts  to be the most important gap in the Enriques-Kodaira classification table for complex surfaces. The standard conjecture concerning this problem states that any minimal class VII surface with $b_2>0$ has $b_2$ curves. By the results of \cite{Ka1} - \cite{Ka3}, \cite{Na1} - \cite{Na3}, \cite{DOT}, \cite{OT}  this conjecture (if true)  would solve the classification problem completely. We explain a new approach (based on techniques from Donaldson theory) to prove existence of curves on class VII surfaces, and we present recent results obtained using this approach. 
\end{abstract}
\maketitle

This survey is a written version of the talk given by the author at the Postnikov Memorial Conference, Bedlewo (Poland), June 2007. 

\section{Introduction}

 By definition, a  K\"ahlerian surface is a compact complex surface admitting a K\"ahler metric.   A fundamental theorem in complex geometry states that a compact complex surface is  K\"ahlerian if and only if $b_1(X)$ is even.  This follows indirectly from Siu's result \cite{Si}  about the K\"ahlerianity of K3 surfaces and standard  results on the classification of complex surfaces.  A direct proof has  been given by Buchdahl (see \cite{Bu2}, \cite{BHPV}). 

According to the Enriques-Kodaira classification table \cite{BHPV}, if $X$ is a non-K\"ahlerian surface, i.e. its first Betti number $b_1(X)$ is odd, then its minimal model $X_{\rm min}$ belongs to one of the following three classes:
\begin{enumerate}
\item Primary and secondary Kodaira surfaces,
\item Non-K\"ahlerian properly elliptic  surfaces, 
\item Minimal class VII surfaces.
\end{enumerate}

A primary Kodaira surface is a topologically non-trivial, locally trivial principal elliptic fiber bundle over an elliptic base. Such a surface has $b_1=3$. A secondary Kodaira surface is a surface with $b_1=1$ which admits a primary Kodaira surface as unramified covering. A properly elliptic surface is an elliptic surface with Kodaira dimension $\kod(X)=1$. Finally, a class VII surface is a complex surface $X$ having $b_1(X)=1$ and $\kod(X)=-\infty$. The condition $\kod(X)=-\infty$ implies $p_g=0$ hence, using standards identities  between analytical and topological invariants of complex surfaces (see Theorem 2.7 p. 139 \cite{BHPV}) one obtains $b_+(X)=0$. Therefore class VII surfaces are interesting from a differential topological point of view: they are 4-manifolds with $b_1=1$ and negative definite intersection form. The simplest examples of class VII surfaces are primary Hopf surfaces: a primary Hopf surface is the quotient of the punctured plane $\C^2\setminus\{0\}$ by an infinite cyclic group $H$ which acts properly discontinuously by holomorphic transformations. For instance, one can take $H=\langle T\rangle$, where $T(z_1,z_2)=(\alpha_1z_1,\alpha_2z_2)$ with $|\alpha_i|<1$. Any primary Hopf surface is diffeomorphic to $S^1\times S^3$. From a historical point of view these are the first examples  of non-K\"ahlerian compact manifolds.     

Note that  the surfaces in the first two classes of the list above are all elliptic, so they are well understood and can be classified using classical  complex geometric methods.
On the other hand, the third class in the list above -- the  mysterious  class VII  -- is not understood yet  and, although much progress has recently been obtained, it resists since many decades to the efforts of the experts. This is considered to be the most important gap in the Enriques-Kodaira classification table.
One can wonder why is the classification problem for class VII surfaces so difficult. Probably the main source of difficulty is the lack of lower dimensional complex geometric sub-objects. Indeed, it is not known (and till recently there was no available method to decide) whether a class VII surface contains a holomorphic curve or a holomorphic foliation.
The classification problem can be understood in two ways, a weak way  and a strong way:
\begin{itemize}
\item [A.] Classify class VII surfaces up to deformation equivalence,
\item [B.] Classify class VII surfaces up to biholomorphic equivalence, and describe explicitly the ``moduli  space"  which corresponds to a given deformation class.

\end{itemize}

Note that such a  ``moduli space" might be a highly non-Hausdorff holomorphic stack. The weaker problem A is   interesting from a topological point of view, because answering this problem would solve the following  differential topological problem: \\

{\it What are the possible diffeomorphism types of compact, connected, oriented 4-manifolds $M$  with $b_1(M)=1$ and negative definite intersection form which admit an integrable almost complex structure? }\\

Interestingly, even the much simpler question ``What are the possible isomorphism types of fundamental groups  $\pi_1(M)$ of such 4-manifolds" is a very difficult, still non-answered question. \\

Although very difficult, there is much hope that in fact the classification problem of class VII surfaces has a simple answer; in other words it is believed that this class is actually quite small. This guess is supported by the {\it known} classification of class VII surfaces with $b_2=0$ (as we will see below)  and  the remark that, in the K\"ahlerian case ($b_1(X)$ even), the condition $\kod(X)=-\infty$  determines the very small class  of surfaces which are rational or ruled.\\

A class VII surface with $b_2=0$ is biholomorphic to either a Hopf surface (i.e. a compact surface whose universal cover is $\C^2\setminus\{0\}$) or an Inoue surface (i.e. a class VII surface which is the free quotient of $\C\times H$ by a properly discontinuous affine action). 
This result was first stated by Bogomolov \cite{Bo1}, \cite{Bo2}, and is usually called Bogomolov's theorem, but  it appears that the first complete proofs have been given in \cite{Te1} and \cite{LYZ}. Both Hopf and Inoue surface are well understood and completely classified.   

There are many examples  of minimal class VII surfaces with $b_2>0$. By results of Kato \cite{Ka1}, \cite{Ka2}, \cite{Ka3} and Dloussky \cite{D1} all these surfaces contain a global spherical shell (a GSS), i.e. an open submanifold $\Sigma$  which is biholomorphic to a standard neighborhood of $S^3$ in $\C^2$ and does not separate $X$ (i.e. $\pi_0(X\setminus U)$ is trivial). Minimal class VII surfaces with $b_2>0$ containing a GSS (also called GSS surfaces, or Kato surfaces) are well understood.   
For instance, an important result of Kato states that  
\begin{thry}
Every GSS surface contains $b_2$ rational curves and is a global deformation (a degeneration) of a 1-parameter family of blown up primary Hopf surfaces. 
\end{thry}

  In particular, all GSS surfaces with fixed second Betti number $b>0$ are deformation equivalent, and are diffeomorphic to   $(S^1\times S^3)\#b \bar\P^2$.

 The biholomorphic classification of GSS surfaces is also well understood, and  certain moduli spaces of such surfaces have recently been described explicitly \cite{OT}.  Therefore the following conjecture formulated by Nakamura \cite{Na2} would solve in principle the classification problem for class VII surfaces:
\\
\\
{\it The GSS conjecture: Every minimal class VII surfaces $X$ with $b_2(X)>0$ has a GSS, so it belongs to the list of known class VII surfaces.}\\

If the conjecture is true, all minimal class VII surfaces $X$ with a fixed second Betti number will be deformation equivalent and diffeomorphic to $(S^1\times S^3)\#b \bar\P^2$ (so not very interesting from a differential topological point of view).   

Since a spherical shell is a non-compact object, it seems that establishing the existence of a GSS is  a very difficult task. Fortunately, a recent result of Dloussky-Oeljeklaus-Toma \cite{DOT} reduces the problem  to the existence of sufficiently many rational curves.
\begin{thry} If a minimal class VII surface with $b_2>0$ has $b_2$ rational curves, then it is a GSS surface.
\end{thry}

Any GSS surface contains exactly $b_2$ rational curves, but the intersection and self-intersection numbers of these curves  are not invariant under deformation. For instance, there exists a holomorphic  family $(X_z)_{z\in D}$ of GSS surfaces with $b_2=1$ such that for $z\ne 0$ the only irreducible curve of  $X_z$  is a homologically trivial singular  rational curve $C_z$, whereas the only irreducible curve of $X_0$ is a singular rational curve $D_z$ with $D_z^2=-1$.  As $z\to 0$ the volume of $C_z$ tends to infinity, so the closed positive current associated with $C_z$ does not converge to a current in $X_0$ as $z\to 0$.  This shows that in our non-K\"ahlerian framework existence of curves cannot be established using Gromov-Witten type invariants.    Gromov-Witten theory cannot be extended to the non-K\"ahlerian complex geometric framework, because on a non-K\"ahler manifold the volume of the curves in a fixed homology class cannot be a priori bounded.

Any GSS surface contains  a {\it cycle of rational curves}, i.e.    an effective divisor $\sum_{i=1}^k C_i$ ($1\leq k\leq b_2$), where either $k=1$ and $C_1$ is a singular rational curve with a normal crossing, or  $k\geq 2$, $C_i$ are smooth rational and 
$$C_1\cdot  C_2=\dots= C_{k-1}\cdot C_k=C_k\cdot C_{1}=1\ ,\ C_i\cdot C_j=0\ {\rm for}\ j-i\ne\pm 1\ {\rm mod}\ k.$$
 In order to solve the weaker classification problem (up to deformation equivalence) it suffices to prove the existence of a {\it cycle of rational curves}. Indeed, by a theorem of Nakamura \cite{Na1}, one has
\begin{thry} A minimal class VII surface which contains a cycle of rational curves is a  global deformation (a degeneration) of a 1-parameter family of blown up primary Hopf surfaces. 
\end{thry}
 
In conclusion, in order to solve the weak classification problem for class VII surfaces it suffices to prove that \\
\\
{\it Conjecture A:  Any minimal class VII surface $X$ with $b_2>0$ has a cycle of rational curves,}\\
\\
whereas  in order to solve the strong classification problem it suffices to prove that \\
\\
{\it Conjecture B: Any minimal class VII surface $X$ with $b_2>0$ has $b_2$ rational curves.}
\\

In the articles \cite{Te2} and \cite{Te3}, we proved  that the existence of curves on class VII surfaces can be established (at least for surfaces with  small second Betti number) using a combination of gauge theoretical and complex geometric techniques.  The main idea is to prove that if a minimal  class VII surface  with $b_2>0$  had no curves, a certain moduli space of polystable rank 2-bundles ($PU(2)$-instantons)  on $X$    would contain a smooth compact  component consisting of stable bundles (irreducible instantons). The existence of such a component leads to a contradiction. This method yields the following result which solves the strong classification problem in the case $b_2=1$ and the weak  classification problem in the case $b_2=2$.
\begin{thry} Conjecture B is true for $b_2=1$. Conjecture A is true for $b_2=2$.
\end{thry}
The purpose of this article is to explain in a geometric, non-technical way the idea of the proof,  pointing out the role of the topological and gauge theoretical techniques.

\section{Instantons on non-K\"ahlerian surfaces}

Let $E$ be a Hermitian rank 2-bundle on a (compact, connected, oriented) Riemannian 4-manifold $(M,g)$. Let $D=\wedge^2E$  be the determinant line bundle of $E$ and $a$ a fixed Hermitian connection of $D$.  We denote by ${\cal A}_a(E)$ the space of $a$-oriented connections, i.e. of Hermitian connections $A$ on $E$ which induce $a$ on the determinant line bundle $D$. Our gauge group is   ${\cal G}:=\Gamma(M,SU(E))$. Let 
$${\cal A}_a^\ASD(E):=\{A\in {\cal A}_a(E)|\ (F_A^0)^+=0\}$$ 
be the  subspace of projectively ASD $a$-oriented connections (or $a$-oriented instantons),  
$${\cal M}_a^\ASD(E):=\qmod{{\cal A}_a^\ASD(E)}{{\cal G}}$$
 the corresponding moduli space, and ${\cal M}_a^\ASD(E)^*\subset {\cal M}_a^\ASD(E)$  the open subspace consisting of  orbits of irreducible projectively ASD oriented connections. 
 
 Denoting by $P$ the $PU(2)$-bundle associated with the unitary  frame bundle of $E$, one has obvious identifications
 $${\cal A}_a(E)={\cal A}(P)\ ,\ {\cal A}_a^\ASD(E)={\cal A}^\ASD(P)\ ,\ {\cal G}=\Gamma(M,P\times_{\Ad} SU(2))\ , 
 $$
 so $a$-oriented instantons can be regarded as $PU(2)$-instantons, and the spaces ${\cal A}_a^\ASD(E)$,  ${\cal M}_a^\ASD(E)$ are independent of $a$ up to canonical isomorphisms. Note however that in general the natural morphism ${\cal G}/\Z_2\to \Aut(P)$ is not surjective, so our moduli space ${\cal M}_a^\ASD(E)$ cannot be identified with the usual  moduli space ${\cal M}^\ASD(P)={\cal A}^\ASD(P)/\Aut(P)$ of $PU(2)$-instantons on $P$. The latter moduli space is the quotient of the former by the natural action of $H^1(M,\Z_2)$ given by tensoring with flat $\{\pm 1\}$-connections. 
 \begin{re}\label{compact} Put $\Delta(E):=4c_2(E)-c_1(E)^2=-p_1(P)$. Then ${\cal M}_a^\ASD(E)$ is
 \begin{enumerate}
 \item  empty  when $\Delta(E)<0$,
 \item compact when $\Delta(E)\leq 3$.
 \end{enumerate}
 \end{re}
 \pf The first statement follows from the well-known identity relating the energy of a  $PU(2)$-instanton to its Pontrjagin number. The second  follows from the first and Uhlenbeck compactness theorem \cite{DK}.
 \qed

 Suppose now that $(X,g)$ is a complex surface endowed with a Gauduchon metric, i.e. a Hermitian metric $g$ whose associated (1,1)-form $\omega_g$ satisfies the equation $\partial\bar\partial\omega_g=0$. Every conformal class of Hermitian metrics contains a Gauduchon metric which is unique up to constant factor (see \cite{Gau}), so Gauduchon metrics exist on any surface. A Gauduchon metric defines a morphism of Lie groups $\deg_g:\Pic(X)\to \R$ given by
 $$\deg_g({\cal L}):=\int_Xc_1({\cal L},h)\wedge\omega_g\ ,
 $$
 where ${\cal L}$ is a holomorphic line bundle on $X$, $h$ is any Hermitian metric on ${\cal L}$, and $c_1({\cal L},h)$ denotes the Chern form of the Chern connection associated with the pair $({\cal L},h)$. For instance, if ${\cal O}(D)$ is the holomorphic  line bundle associated with an effective divisor $D\subset X$, then $\deg_g({\cal O}(D))=\Vol_g(D)$.

 If $b_1(X)$ is  odd, the degree map is not a topological invariant (i.e. it does not vanish on the connected component $\Pic^0(X)\simeq H^1(X,{\cal O})/H^1(X,\Z)$ of the zero element in $\Pic(X)$). For instance, if $X$ is a class VII surface, then one has a canonical isomorphism $\Pic^0(X)\simeq\C^*$, and the degree of the line bundle ${\cal L}_z$ associated with $z\in \C^*$ is given by $\deg({\cal L}_z)=C_g\log|z|$  for a positive constant   $C_g$  determined by the metric $g$. 
 
For a coherent sheaf ${\cal F}$ on $X$ one defines $\deg_g({\cal F}):=\deg_g(\det({\cal F}))$. A holomorphic rank 2-bundle ${\cal E}$ on $X$ is called {\it stable} if $\deg_g({\cal L})<\frac{1}{2}\deg({\cal E})$ for every sheaf monomorphism ${\cal L}\to {\cal E}$ with torsion free quotient,   and {\it polystable} if it is either stable or isomorphic to  the direct sum of  two line bundles of the same degree. 

Fix now a ${\cal C}^\infty$ Hermitian rank 2-bundle $(E,h)$ on $X$, and fix a holomorphic structure ${\cal D}$ on its determinant  line bundle $D$. We denote by ${\cal M}^\st_{\cal D}(E)$, ${\cal M}^\pst_{\cal D}(E)$ the moduli sets of stable (respectively   polystable) holomorphic structures ${\cal E}$ on $E$ which induce ${\cal D}$ on the determinant line bundle (modulo the equivalence relation defined by the ${\cal G}^\C:=\Gamma(X,SL(E))$-action on the space of these structures).  ${\cal M}^\st_{\cal D}(E)$ is a Hausdorff open subspace of  the moduli space ${\cal M}^\s_{\cal D}(E)$ of {\it simple} holomorphic structures ${\cal E}$ on $E$ with $\det({\cal E})={\cal D}$, which has a natural structure of a finite dimensional (in general  non-Hausdorff) complex space \cite{LO}, \cite{LT1}. The Kobayashi-Hitchin correspondence for complex surfaces \cite{Do1}, \cite{Bu1}, \cite{LY}, \cite{LT1}, \cite{LT2} states that 
\begin{thry}
Let $a$ be the Chern connection associated with the holomorphic structure ${\cal D}$ and the Hermitian metric $\det(h)$ on $D$. The assignment  $A\mapsto \bar\partial_A$ defines a real analytic isomorphism ${\cal M}_a^\ASD(E)^*\stackrel{\simeq}{\to}{\cal M}^\st_{\cal D}(E)$ which extends to a bijection ${\cal M}_a^\ASD(E)\stackrel{\simeq}{\to}{\cal M}^\pst_{\cal D}(E)$.
 \end{thry}
  
  We will endow ${\cal M}^\pst_{\cal D}(E)$ with the topology which makes this bijection a homeomorphism. This important theorem furnishes a tool to compute moduli spaces of instantons on complex surfaces using complex geometric methods.  The point here is that classifying  (poly)stable holomorphic  bundles is much easier than classifying the solutions of the ASD equation (which is a difficult non-linear PDE system). This principle plays a fundamental role in Donaldson  theory (see \cite{Do1}, \cite{DK}) and made possible the first explicit computations of Donaldson invariants.
Unfortunately  classifying holomorphic bundles  on {\it non-algebraic} surfaces  becomes a very difficult (sometimes hopeless) problem. The difficulty (specific to the non-algebraic framework) is the  appearance of non-filtrable holomorphic bundles \cite{BLP}.  A holomorphic rank 2-bundle ${\cal E}$ is called {\it filtrable} if its associated locally free sheaf (denoted by the same symbol) contains a coherent  rank 1-subsheaf  or, equivalently, if it  fits in an exact sequence of the form
$$0\map {\cal L}\map {\cal E}\map {\cal M}\otimes{\cal J}_Z\map 0\ ,
$$
where ${\cal L}$, ${\cal M}$ are holomorphic line bundles and $Z\subset X$ is a zero-dimensional locally complete intersection.  On algebraic surfaces all bundles are filtrable, so classifying holomorphic rank 2-bundles reduces to classifying locally free extensions of    ${\cal M}\otimes{\cal J}_Z$ by ${\cal L}$, where ${\cal L}$, ${\cal M}$ vary in $\Pic(X)$, and $Z$ in the Douady space of   zero-dimensional locally complete intersections.

\section{Existence of curves on class VII surfaces}

Let $X$ be a minimal class VII surface with $b_2(X)=b>0$. We get easily $c_2(X)=e(X)=b$ and (using Noether's formula  \cite{BHPV})  $c_1(X)^2=-b$. We denote by ${\cal K}$ the canonical holomorphic line bundle on $X$, and by $K$ its underlying differentiable line bundle. In order to avoid technical difficulties, we make the following {\it simplifying assumptions}:
\begin{enumerate}
\item $\pi_1(X)\simeq\Z$,
\item There exists a Gauduchon metric $g$ on $X$ such that $\deg_g({\cal K})<0$.
\end{enumerate}

The first assumption will simplify our arguments in the following way: First it guarantees that $H^2(X,\Z)$ is a free $\Z$-module.  By Donaldson first theorem \cite{Do2}, there exists a basis $(e_1,\dots,e_b)$ in $H^2(X,\Z)\simeq\Z^b$ which is orthonormal with respect to the intersection form $q_X$, i.e. $e_i\cdot  e_j=-\delta_{ij}$. Taking into account that $c_1(K)^2=-b$ and that $c_1(K)$ is a characteristic class, it is easy to see that (replacing some $e_i$'s by their opposite if necessary) this  basis can be chosen such that $\sum_i e_i=c_1(K)$. Second, the first assumption implies that there exists a unique non-trivial representation $\rho:\pi_1(X,x_0)\to\{\pm 1\}$, so $X$ has a unique double cover $\pi_\rho:\tilde X\to X$. 

The second condition is not restrictive: one can prove that if $\deg_g({\cal K})>0$ for every Gauduchon metric $g$ on $X$,  then the Chern class $c_1({\cal K})_{\rm BC}$ in Bott-Chern cohomology will be pseudo-effective; this implies that $X$ has a numerically pluri-canonical divisor, which is not possible on minimal class VII surfaces with $b_2>0$, by a result of Nakamura (see  \cite{Te3} for details).

Let $(E,h)$ be a Hermitian rank 2-bundle on $X$ with $c_2(E)=0$, $\det(E)=K$. We are interested in the moduli space
$$ {\cal M}:={\cal M}^\pst_{\cal K}(E)\simeq {\cal M}_a^\ASD(E)\ ,$$
where $a$ is the Chern connection associated with the pair $({\cal K},\det(h))$. Since $c_1(E)=-c_1(X)$ we get  $\Delta(E)=b$, so the complex expected dimension of our moduli space is $b$. Moreover, by Remark \ref{compact} 
\begin{re} ${\cal M}$ is compact when $b\leq 3$. 
\end{re}

This result is an easy consequence of the Kobayashi-Hitchin correspondence and Uhlenbeck compactness theorem, but it cannot be obtained using complex geometric arguments (although polystable bundles are defined in complex geometric terms).

The idea    to consider this moduli space for our purposes might look surprising, because there is no obvious relation  between instantons and holomorphic curves.  In order to clarify in a non-technical way how this moduli space will be used, let us first classify the {\it filtrable} rank 2-bundles ${\cal E}$ with $c_2({\cal E})=0$, $\det({\cal E})={\cal K}$. Such a bundle fits in an exact sequence $0\to{\cal L}\to{\cal E}\to {\cal K}\otimes{\cal L}^{-1}\otimes {\cal J}_Z\to 0$.
Write $c_1({\cal L})=\sum n_i e_i$. The relation $c_2({\cal E})=0$ gives $|Z|+\sum_i n_i(n_i-1)=0$, which shows that $Z=\emptyset$ and $c_1({\cal L})=e_I:=\sum_{i\in I} e_i$ for a subset $I\subset I_0:=\{1,\dots,b\}$. Therefore for every such index set $I\subset I_0$ we get a class of extensions of the form
\begin{equation}\label{ext}
0\map {\cal L}\textmap{j}{\cal E}\textmap{q} {\cal K}\otimes{\cal L}^{-1}\map 0
\end{equation}
with $c_1({\cal L})=e_I$. We will call them extensions of type $I$. We have proved that any filtrable rank 2-bundle with $c_2({\cal E})=0$, $\det({\cal E})={\cal K}$ is the central term of such an extension. Unfortunately, in general, the same filtrable bundle can be written as an extension in many different ways, so it is not clear at all which bundles obtained in this way are stable. 
An important example is obtained taking $I=I_0$ and ${\cal L}={\cal K}$. For this choice there exists an (essentially unique) non-trivial extension
\begin{equation}\label{A}
0\map {\cal K}\textmap{i} {\cal A}\textmap{p} {\cal O}\map 0\ ,
\end{equation}
because ${\rm Ext}^1({\cal O},{\cal K})=H^1({\cal K})\simeq H^1({\cal O})^\vee\simeq \C$, by Serre duality and Riemann-Roch theorems. If the central term ${\cal A}$ can be written as an extension in a different way, we get a bundle monomorphism $v:{\cal L}\to {\cal A}$ whose image  is not contained in ${\cal K}$. Thus $v$ defines a sheaf monomorphism $p\circ  v:{\cal L}\to {\cal O}$ which cannot be an isomorphism, because, if it were, $v$ would define a splitting of (\ref{A}). Therefore the image $(p\circ v)({\cal L})$ can be identified with the ideal sheaf ${\cal O}(-C)$ of a {\it non-empty} effective divisor $C$. A careful examination of the possible cases which can occur \cite{Te3} shows that $C$   contains either an elliptic curve (in which case $X$ is a GSS surface by Nakamura's results)  or a cycle of rational curves. Therefore
\begin{thry}\label{stableA} If $X$ does not contain a cycle of rational curves, the bundle ${\cal A}$ (defined as the essentially unique non-trivial extension of ${\cal O}$ by ${\cal K}$) cannot be written as an extension in a different way. If, moreover, $\deg_g({\cal K})<0$, the unique sub-line bundle ${\cal K}$ of ${\cal A}$ does not destabilize ${\cal A}$, so ${\cal A}$ is stable.
\end{thry}

It is easy to see that (under our assumptions) there are two non-trivial extensions of type $I_0$: ${\cal A}$, and ${\cal A}':={\cal A}\otimes {\cal L}_\rho$, where ${\cal L}_\rho$ is the flat line bundle defined by the representation $\rho$. The line bundle ${\cal L}_\rho$ is a non-trivial square root of the trivial line bundle ${\cal O}$.\\

We have seen that an easy way to construct holomorphic  rank 2-bundles is to look at extensions of line bundles. But there is another simple, classical  method: taking the push-forward of a line bundle ${\cal N}$ defined on a double cover $\pi:\tilde X\to X$.  One has an obvious isomorphism  $\pi^*(\pi_*({\cal N}))\simeq{\cal N}\oplus \iota^*({\cal N})$, where $\iota:\tilde X\to\tilde X$ denotes the canonical involution of the double cover $\tilde X$. Note that the two summands ${\cal N}$, $\iota^*({\cal N})$
 have the same degree with respect to a pull-back Gauduchon metric $\pi^*(g)$, so $\deg_g(\pi_*({\cal N}))=\deg_{\pi^*(g)}({\cal N})$.  It follows easily that the push-forward bundle  $\pi_*({\cal N})$  is  polystable with respect to any Gauduchon metric $g$ on $X$, and is stable when  ${\cal N}\not\in\pi^*(\Pic(X))$.    

The irreducible  instanton $A$ corresponding to a stable bundle $\pi_*({\cal N})$ obtained in this way will be irreducible, but its pull-back to $\tilde X$ is reducible. Its orbit $[A]\in {\cal M}_a^\ASD(E)$ is a fixed point of the involution defined  by tensoring with the flat $\{\pm\}$-connection associated with the double cover $\pi$. Such instantons are called sometimes {\it twisted reductions} (see \cite{KM}).
Note that in many cases the stable bundles obtained in this way are non-filtrable.

In our case, by our first simplifying assumption, we have a  unique (up to isomorphism) double cover $\pi_\rho:\tilde X\to X$, and one can easily classify the line bundles ${\cal N}\in\Pic(\tilde X)$ for which $c_2(\pi_*({\cal N}))=0$, $\det(\pi_*({\cal N}))\simeq{\cal K}$. One finds in this way $2^{b-1}$ isomorphism classes of stable bundles (irreducible instantons), which are precisely the fixed points of the involution $\otimes{\cal L}_\rho$ acting on our moduli space ${\cal M}$. A natural question is whether these bundles are filtrable. Equivalently, one can ask if, for a given index set $I\subset I_0$, the (isomorphism class of the) central term of an extension of type $I$ can be  a fixed point of the involution $\otimes{\cal L}_\rho$. For instance, if   ${\cal A}$ is such a fixed point, one gets a   diagram
$$
\begin{array}{ccccccccc}
0&\map&{\cal K}&\textmap{i}&{\cal A}&\textmap{p}&{\cal O}&\map&0\\
&&&&\simeq\downarrow\psi\\
0&\map&{\cal K}\otimes{\cal L}_\rho&\textmap{i'}&{\cal A}\otimes{\cal L}_\rho&\textmap{p'}&{\cal L}_\rho&\map&0
\end{array}
$$
with $\psi$ an isomorphism. If $p'\circ \psi\circ i$ vanished, it would follow that $\psi({\cal K})\subset {\cal K}\otimes{\cal L}_\rho$; but $\Hom({\cal K},{\cal K}\otimes{\cal L}_\rho)=H^0({\cal L}_\rho)=0$, because ${\cal L}_\rho$ is a non-trivial holomorphic line bundle of degree 0.  Therefore, one would get $\resto{\psi}{{\cal K}}=0$, which is not possible, because $\psi$ is an isomorphism. This means that $p'\circ \psi\circ i\ne 0$, so one gets a nontrivial section $s\in H^0({\cal K}^\vee\otimes{\cal L}_\rho)$, whose vanishing locus would be a numerically anticanonical divisor and $X$ will be a GSS surface by the main result in \cite{D2}. Therefore
\begin{pr}\label{fixedA}
If ${\cal A}$ is a fixed point of the involution $\otimes{\cal L}_\rho$ then $X$ is a GSS surface.
\end{pr}

\subsection{The case $b=1$}

In the case $b=1$, we get a compact  moduli space ${\cal M}$ of complex dimension 1, which contains the following remarkable subspaces:
\begin{enumerate}
\item The subspace $\Rg$ of reductions.

The reductions in ${\cal M}$ correspond to split polystable bundles of the form ${\cal L}\oplus({\cal K}\oplus{\cal L}^\vee)$ with ${\cal L}\in\Pic^0(X)\simeq\C^*$ and $\deg_g({\cal L})=\frac{1}{2}\deg_g({\cal K})$. Taking into account the explicit form of the degree map $\resto{\deg_g}{\Pic^0(X)}$ explained above, we see that $\Rg$ is a circle.
\\
\item The subspace of stable type $\emptyset$-extensions, i.e. of stable extensions of the form 
$$0\map{\cal L}\map {\cal E}\map {\cal K}\otimes{\cal L}^{-1}\map 0\ ,
$$
with ${\cal L}\in\Pic^0(X)$.  The stability condition  implies $\deg_g({\cal L})<\frac{1}{2}\deg_g({\cal K})$, and one can prove easily as in the proof of Theorem \ref{stableA}  that,   for every ${\cal L}\in\Pic^0(X)$ with $\deg_g({\cal L})<\frac{1}{2}\deg_g({\cal K})$ there exists an essentially unique extension ${\cal E}_{\cal L}$ of this type, which is stable. Moreover, using that $X$ is minimal, one can prove that the assignment ${\cal L}\mapsto{\cal E}_{\cal L}$ is injective.  The subspace of these stable bundles can be identified with a punctured disk $D^\bullet$.
\begin{figure}[h]
\centering
\scalebox{0.4}
{\includegraphics{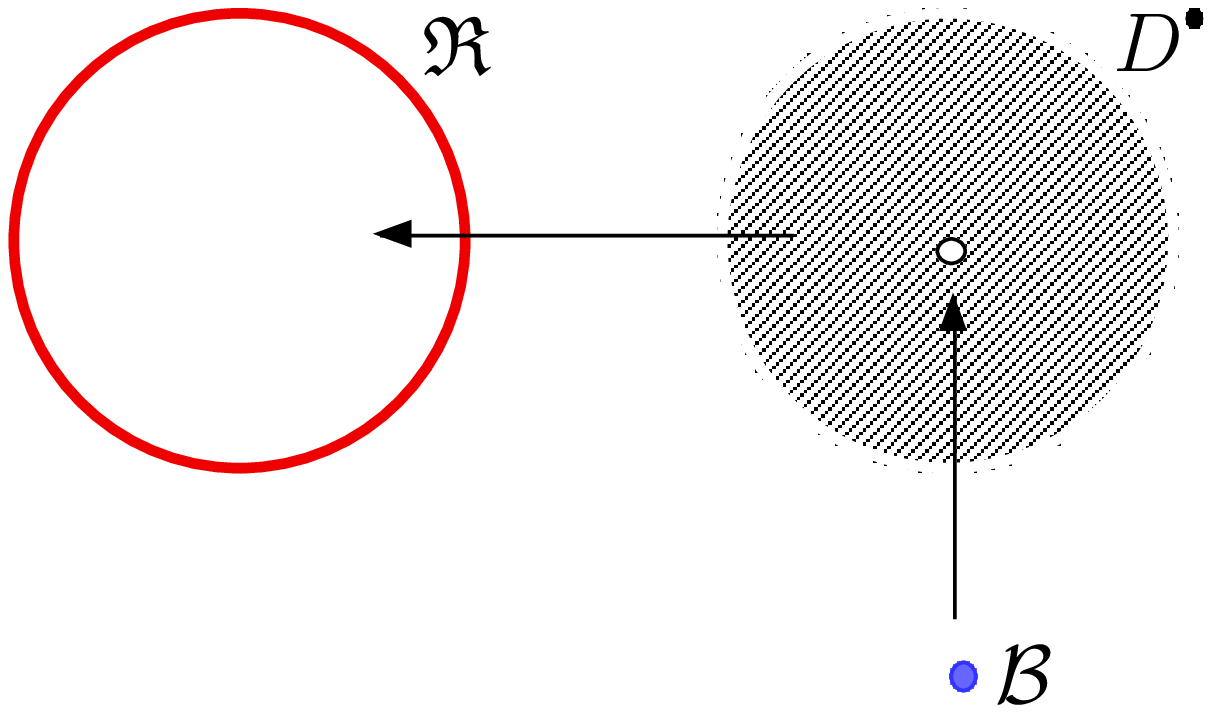}}
\label{picture1}
\end{figure}
\item A push-forward  bundle  ${\cal B}=[\pi_\rho]_*({\cal N})$, which is the unique fixed point of the involution $\otimes{\cal L}_\rho$.
 \end{enumerate}
%


 We   know that ${\cal M}$ must be compact and that ${\cal M}^\st={\cal M}\setminus\Rg$ is a 1-dimensional complex space. Using a well-known smoothness criterion for the universal deformation of a simple bundle,  it follows easily that the complex space ${\cal M}^\st$ is   smooth. On the other hand the general deformation theory for instantons \cite{DK} shows that  ${\cal M}$ has the structure of a Riemann surface with boundary at the points of $\Rg$.
 
 One can ask now what are the incidence relations between the three subspaces listed above. It is easy to see that the circle  $\Rg$ compactifies $D^\bullet$ towards its outer end, whereas $\{{\cal B}\}$ compactifies $D^\bullet$ towards the origin. Indeed, the only way to compactify a closed punctured disk to get a smooth Riemann surface with boundary is to add a point at the place of the missing origin. This point must be a fixed point of $\otimes{\cal L}_\rho$ by Browder fixed point theorem (because   $\otimes{\cal L}_\rho$  acts freely on $D^\bullet\cup\Rg$),  so it must coincide with ${\cal B}$.  Therefore the first three pieces of our moduli space fit together in the obvious way and form a closed disk ${\cal M}_0$ bounded by the circle of reductions $\Rg$. ${\cal M}_0$ is connected component of ${\cal M}$. \vspace{-2mm}
 \begin{figure}[h]
\centering
\scalebox{0.4}
{\includegraphics{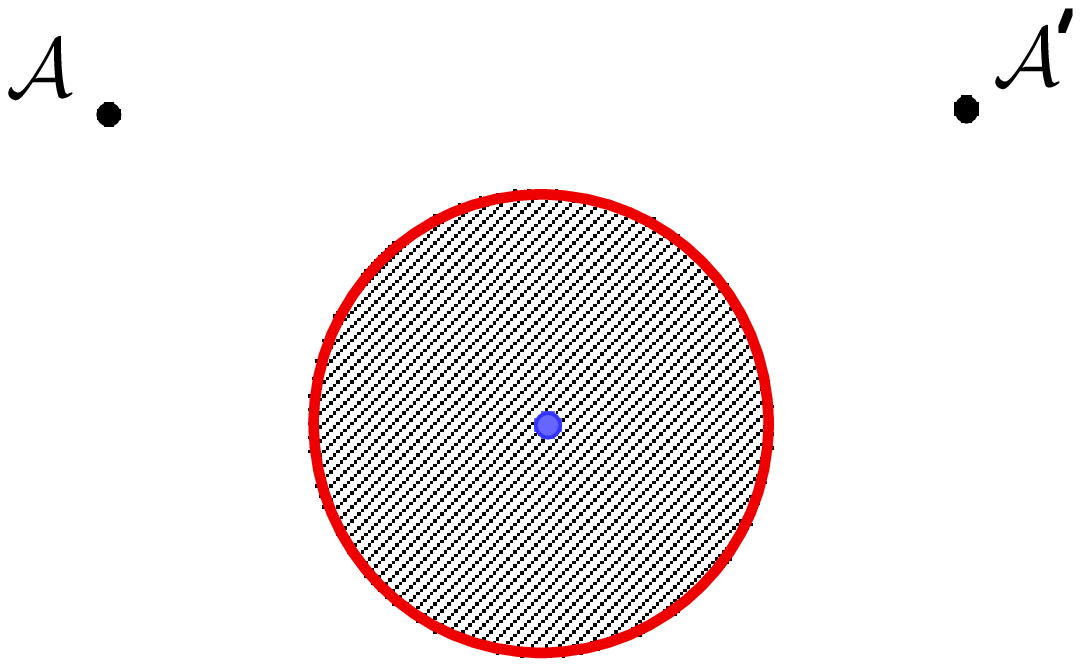}}
\label{picture2}
\end{figure} 
\vspace{-2mm}
 
 We want to prove that $X$ has a cycle. If not, by Theorem \ref{stableA}, we get two new points in ${\cal M}$, namely the stable filtrable bundles ${\cal A}$ and ${\cal A}'$.  ${\cal A}$, ${\cal A}'$ and ${\cal E}_{\cal L}$ ($\deg_g({\cal L})<\frac{1}{2}\deg_g({\cal K})$) are  the only filtrable stable bundles in our moduli space. I came to the idea to use instantons (in order to prove new properties of class VII surfaces) as I was trying to understand how these points fit in the moduli space ${\cal M}$.

 Obviously, one of the following must hold:
 \begin{enumerate}
 \item ${\cal A}\in D^\bullet\cup\{{\cal B}\} $,
 \item  The connected component of ${\cal A}$ is a closed Riemann surface $Y\subset {\cal M}^\st$, and all the points of $Y\setminus\{{\cal A},{\cal A}'\}$ correspond to non-filtrable bundles.
\end{enumerate} 

Note that ${\cal K}\not\in \Pic^0(X)$, so  ${\cal A}\in D^\bullet$ would imply that the bundle ${\cal A}$ has another representation as an extension. Therefore the first case implies the existence of a cycle by Theorem \ref{stableA} and Proposition \ref{fixedA}. On the other hand the second case leads to a contradiction by the following
\begin{thry}  Let $X$ be a surface with $a(X)=0$, $E$ a differentiable rank 2-bundle on $X$, ${\cal D}$ a holomorphic structure on $\det(E)$, and $Y$ a Riemann surface. There exists no holomorphic map $\psi:Y\to {\cal M}^\s_{\cal D}(E)$ whose image contains both filtrable and non-filtrable bundles.  
\end{thry}

\pf The first step in the proof of this result is to remark that a holomorphic map $\psi:Y\to {\cal M}^\s(E)$ defined on a Riemann surface is induced by a holomorphic family $({\cal F}_y)_{y\in Y}$ of rank 2-bundles on $X$ parameterized by $Y$, i.e. a {\it universal} rank 2-bundle ${\cal F}$ on $Y\times X$. In general there might be an obstruction in $H^2({\cal O}_Y^*)$ to the existence of such a universal bundle, but this group vanishes for a Riemann surface $Y$ (see \cite{Te2} for details).   
Next, one uses the fact that $X$ has  no non-constant meromorphic functions (so it is ``very" non-algebraic), whereas $Y$ is a compact, complex curve, so it is algebraic.  This ``incompatibility" between the  two factors  implies that  very few holomorphic bundles on their product can exist. In order to see this, we will switch the roles of the two factors, namely we regard  ${\cal F}$ as a family $({\cal F}^x)_{x\in X}$ of bundles on $Y$ parameterized  by $X$. 

We explain now the proof of the theorem in a particular case which illustrates very well  how this ``incompatibility" between the algebraicity properties of the two factors is used. Suppose that ${\cal F}^x$ is stable on $Y$ for generic $x\in X$. In this case one gets a meromorphic map $u:X\dashrightarrow {\cal M}^\sst$ in a moduli space of semistable bundles over $Y$, and this moduli space is a projective  variety. Therefore, since $a(X)=0$, this meromorphic map must be constant, so there exists a stable bundle ${\cal F}^0$ on $Y$ such that ${\cal F}^x\simeq {\cal F}^0$ for every $x$ in a  non-empty  Zariski open subset   $X_0\subset X$. The sheaf ${\cal T}:=(p_X)_*(p_Y^*({\cal F}^0)^\vee\otimes{\cal F})$ is a rank 1-sheaf on $X$ and comes with a tautological morphism $\tau:p_X^*({\cal T})\otimes p_Y^*({\cal F}_0)\to {\cal F}$ which is an isomorphism on $Y\times X_0$. The restriction to a fiber $\{y\}\times X\simeq X$ yields a morphism ${\cal T}\otimes {\cal F}^0(y)\simeq {\cal T}\oplus{\cal T}\to {\cal F}_y$ whose restriction to  $X_0$ is an isomorphism. Thus all bundles ${\cal F}_y$ admit rank 1-subsheaves, so they are filtrable.
\qed

\subsection{The case $b=2$}

The moduli space ${\cal M}$ is again compact, but this time its complex dimension is 2.  We use a similar strategy, which begins with the description of the geometry of  the component ${\cal M}_0$ of the moduli space ${\cal M}$ consisting of reductions, stable filtrable bundles defined by extensions of type $I{\subset}_{\hskip-6pt\raisebox{-3.5pt}{${\scriptscriptstyle \ne}$}} I_0$, and  twisted reductions. This is the ``known"  component of the moduli space. Under our simplifying assumptions we have two circles  or reductions $\Rg'$, $\Rg''$, a 2-dimensional family of extensions of type $\emptyset$, two 1-dimensional families of extensions of type $\{e_1\}$ and $\{e_2\}$, and two  twisted  reductions  ${\cal B}_1$, ${\cal B}_2$. However this time, it is much more difficult to understand (and prove!) how these strata fit together \cite{Te3}, and to see why their union is a connected component of the moduli space, as claimed.  

As in the case $b_2=1$ we ask whether the type $I_0$ extension ${\cal A}$  belongs to the known component ${\cal M}_0$  component or not. If it does, we get a cycle in $X$, by Theorem \ref{stableA} and Proposition \ref{fixedA}. If not, we get a smooth compact surface $Y\subset {\cal M}^\st$ such that all the points of $Y\setminus\{{\cal A},{\cal A}'\}$  correspond to non-filtrable bundles. Unfortunately this time we cannot get a contradiction as in the case $b_2=1$, because $Y$ might be  non-algebraic.  The proof uses a long and difficult  analysis of the properties of the surface $Y$ (e.g. its intersection form, Chern classes) in the case when, by reductio ad absurdum, ${\cal A}$ did not belong to  ${\cal M}_0$. Using the classification of complex surfaces, one eliminates one by one all the possibilities \cite{Te3}.

\end{document}